\documentclass[10pt,lettersize,reqno,fullpage]{amsart}
\usepackage{amsmath}
\usepackage{dsfont}
\usepackage{bbm}
\usepackage{ulem}
\usepackage{float}
\usepackage{pgfplots}
\usepackage{amsmath,amsfonts,amssymb}   
\usepackage{float}
\usepackage{amsmath,amssymb,amsfonts,mathrsfs,verbatim,enumitem,pstricks,amsthm, graphicx}
\usepackage[all]{xy}
\usepackage{centernot, xr, accents}
\usepackage{hyperref}

\usepackage{hyperref}
\hypersetup{
	colorlinks,
	citecolor=red,
	filecolor=black,
	linkcolor=blue,
	urlcolor=black
}

\providecommand \hf{\hspace*{0.5cm}}
 
\newtheorem{thm}{Theorem}[section]

\newtheorem{lmm}[thm]{Lemma}

\newtheorem{rem}[thm]{Remark}

\theoremstyle{definition}
\newtheorem{defn}[thm]{Definition}

\usepackage{amsmath,stackengine,scalerel}

\makeatletter
\renewcommand{\tocsection}[3]{%
  \indentlabel{\@ifnotempty{#2}{\bfseries\ignorespaces#1 #2\quad}}\bfseries#3}
\renewcommand{\tocsubsection}[3]{%
  \indentlabel{\@ifnotempty{#2}{\ignorespaces#1 #2\quad}}#3}

\newcommand\@dotsep{4.5}
\def\@tocline#1#2#3#4#5#6#7{\relax
  \ifnum #1>\c@tocdepth 
  \else
    \par \addpenalty\@secpenalty\addvspace{#2}%
    \begingroup \hyphenpenalty\@M
    \@ifempty{#4}{%
      \@tempdima\csname r@tocindent\number#1\endcsname\relax
    }{%
      \@tempdima#4\relax
    }%
    \parindent\z@ \leftskip#3\relax \advance\leftskip\@tempdima\relax
    \rightskip\@pnumwidth plus1em \parfillskip-\@pnumwidth
    #5\leavevmode\hskip-\@tempdima{#6}\nobreak
    \leaders\hbox{$\m@th\mkern \@dotsep mu\hbox{.}\mkern \@dotsep mu$}\hfill
    \nobreak
    \hbox to\@pnumwidth{\@tocpagenum{\ifnum#1=1\bfseries\fi#7}}\par
    \nobreak
    \endgroup
  \fi}
\AtBeginDocument{%
\expandafter\renewcommand\csname r@tocindent0\endcsname{0pt}
}
\def\l@subsection{\@tocline{2}{0pt}{2.5pc}{5pc}{}}
\makeatother

\begin{document}

\title[Kac's Central Limit Theorem by Stein's Method]{Kac's Central Limit Theorem by Stein's Method}
\thanks{Declarations of interest: none}

\author[Suprio Bhar]{Suprio Bhar*}
\address{Department of Mathematics and Statistics, Indian Institute of Technology, Kanpur, Uttar Pradesh- 208016, India.}
\email{suprio@iitk.ac.in }
\thanks{*Corresponding Author}

\author[Ritwik Mukherjee]{Ritwik Mukherjee}
\address{School of Mathematical Sciences, National Institute of Science Education and Research, HBNI, Bhubaneswar, Odisha- 752050, India.}
\email{ritwikm@niser.ac.in}

\author[Prathmesh Patil]{Prathmesh Patil}
\address{Universit\'{e} de Gen\`{e}ve, Section de Math\'{e}matiques, UNI DUFOUR,
24, rue du Général Dufour,
Case postale 64,
1211 Geneva 4, Switzerland.}
\email{Prathmesh.Patil@etu.unige.ch}



\date{}

\begin{abstract}
In $1946$, Mark Kac proved a Central Limit type theorem  
for a 
sequence of random variables that were not independent. The random variables under consideration were obtained from the 
angle-doubling map. 
The idea behind Kac's proof was to show that although the 
random variables under consideration were not independent, they were what he calls \textit{statistically independent} 
(in modern terminology, this concept is called long range independence). The final conclusion of his paper was that the sample averages of the 
random variables, suitably normalized converges 
to the standard normal distribution. We describe a new proof of Mark Kac's result by applying Stein's method and show that the normalized sample averages converge to the standard normal distribution in the 
Wasserstein metric, which is stronger than the convergence in distribution.
\end{abstract}
\maketitle

\section{Introduction}
Two fundamental results in Probability Theory are the 
Strong Law of Large Numbers and the Central Limit Theorem (see \cite{Billingsley-book,Chung-book}). 
The former states that the sample averages of a sequence of i.i.d random variables having finite mean, converge almost surely to the mean. The latter states that if in addition, the variance is finite and non zero, then the error is of the order of $\frac{1}{\sqrt{n}}$, in the sense of convergence in distribution, where $n$ is the sample size.\\ 
\hf \hf It is a natural question to ask what is the situation for random variables that are possibly dependent. 
A simple example of that is as follows: choose $\alpha$ in the closed interval $[0,1]$ with the uniform 
distribution. Now consider the sequence of numbers $\{2^n \alpha\}$ modulo $1$. 
Given a reasonable function $f$ on $[0,1]$, we can consider the sequence of random variables $\{f(2^n \cdot)\}$. The Ergodic Theorem states that for an $L^1$ function $f$, almost surely, 
the sample averages converge to the integral of $f$.\\  
\hf \hf In his paper \cite{Kac_Ber}, Mark Kac shows that the rate of this convergence is of the order of 
$\frac{1}{\sqrt{n}}$. 
More precisely, he shows that a suitable normalized sample average converges (in the sense of convergence in distribution) to the 
standard normal distribution. \\ 
\hf \hf We give an alternative proof of Kac's Theorem by using Stein's method.
We show that the normalized sample averages converge to the standard normal distribution in the 
Wasserstein metric, which is stronger than the convergence in distribution.\\ 
\hf \hf Kac's Theorem and main result (Theorem \ref{mt_bmp}) of this article have been discussed in Section \ref{S2:main-results}. In Section \ref{S3:proofs}, we give a detailed proof of Theorem \ref{mt_bmp}.

\section{Kac's Central Limit Theorem and Main results}\label{S2:main-results}
\hf \hf Let us now recall the setup and result of Kac's paper \cite{Kac_Ber}. Consider the ergodic transformation $T:[0,1]\longrightarrow [0,1]$ given by \begin{align}
T(\alpha)&:= 2 \alpha ~~\textnormal{mod} ~~1. \label{Bernouli_shift} 
\end{align}
We shall refer to this map as the `angle-doubling' map (\cite[Example 2.4]{EW}, \cite[Proposition 2.15, pp. 25]{EW}). In the Ergodic Theory literature, the map $T$ is also linked to the Bernoulli Shifts (see \cite{EW}, \cite{KrengelBk} and \cite{ZimmerBk}). 

\begin{center}
\begin{tikzpicture}
\begin{axis}
[
    xlabel={$x$-Graph of map $T$}, 
    ylabel={$T(x)$}, 
    grid = both, 
    xmin=0, 
    xmax=1, 
    ymin=0, 
    ymax=1,
]
\addplot[domain = 0: 0.5, blue, no marks, domain=0:1, samples=50] {2*x};
\addplot[domain = 0.5: 1, blue, no marks, domain=0:1, samples=50] {(2*x)-1};
\draw[dashed, red, line width=0.95pt] ({axis cs: 0.5,0}|-{rel axis cs:0,0}) -- ({axis cs: 0.5,0}|-{rel axis cs:0,1});
\end{axis}
\end{tikzpicture}
\end{center}

\hf \hf Choose a real number $\alpha \in [0,1]$ randomly (with respect to the uniform distribution on $[0,1]$) and consider the iterates $\alpha$, $T \alpha$, $T^2 \alpha, \cdots$. Let $f:[0,1]\longrightarrow \mathbb{R}$ be any $L^1$ function. Consider the random variables 
$X^f_k:[0,1]\longrightarrow \mathbb{R}$, given by 
\begin{align}
X^f_k(\alpha)&:= f(T^k(\alpha)), k = 0, 1, 2, \cdots. \label{Xf_defn}
\end{align}
Now look at the sample average of $f$ composed with iterates of $T$, namely 
\begin{align}
A_n(\alpha)&:= \frac{X^f_0(\alpha)+ X^f_1(\alpha) + \ldots + X^f_{n-1}(\alpha)}{n}. \label{afn_defn}
\end{align}
Since $T$ is Ergodic, it is a consequence of Birkhoff's Ergodic Theorem (\cite[Theorem 2.30, pp. 44]{EW}) that for almost all $\alpha$ in $[0,1]$, the sample average $A_n(\alpha)$ converges to 
$\mu$, where 
\begin{align*}
\mu&:= \int\limits_{0}^{1} f(t) dt. 
\end{align*} 
\hf \hf In his paper \cite{Kac_Ber}, Kac shows that even though these random variables $X^f_k, k = 0, 1, 2, \cdots$ are not independent, a conclusion similar to that of the 
Central Limit Theorem holds (provided we make certain assumptions on the function $f$, in terms of decay of the Fourier coefficients). \\
Let $\sigma_n^2$ denote the variance of the random variable $X^f_0+X^f_1+\ldots + X^f_{n-1}$. In fact,
\begin{align}
\sigma_n^2 &:= \int\limits_{0}^{1}|A_n(t)|^2 dt. \label{sigm_n_sq_integral}
\end{align}
Let us now assume that 
\begin{align}
\lim_{n\longrightarrow \infty}\frac{\sigma_n^2}{n} & = \sigma^2 >0. \label{sigma_n_lim_assumption}
\end{align}
Note that the existence of $\sigma^2$ and positivity are part of the assumption. In particular, $f$ is not identically zero.\\ 
Next, assume the following condition on the Fourier series expansion of $f$, namely 
\begin{align*}
f(t) &\sim \sum_{n=1}^{\infty} a_n \cos(2\pi n t). 
\end{align*}
Note that we extend $f$ periodically on the whole of $\mathbb{R}$; hence $f$ is an even function. As a result the sine terms 
are absent. 
Moreover, assume that
\begin{align}
|a_n|&< \frac{M}{n^{\beta}} \label{Fourier_condition_f}  
\end{align}
for some constant $M$ and for some $\beta>\frac{1}{2}$. \\ 
Let us now define the random variable $W^f_n$ given by 
\begin{align}
W^f_n(\alpha)&:= \frac{A_n(\alpha)-\mu}{\sigma_n/\sqrt{n}}. \label{Wnf_defn}
\end{align}
The main result of Kac's paper is as follows (which is \cite[Theorem 1, pp. 41]{Kac_Ber}):
\begin{thm}
\label{mt_kac}
Let $f:[0,1]\longrightarrow \mathbb{R}$ be an $L^1$ function, with Fourier coefficients satisfying \eqref{Fourier_condition_f}. 
Let $\{X^f_k\}_{k=0}^{\infty}$ be the collection of random variables as defined by \eqref{Xf_defn} and let $W^f_n$ be as defined 
by \eqref{Wnf_defn}. Suppose equation \eqref{sigma_n_lim_assumption} is satisfied. Under these assumptions, $\forall a \in \mathbb{R}$   
\begin{align*}
    \lim_{n\rightarrow\infty}\mathbb{P}( W_n^f\leq a) = \frac{1}{\sqrt{2\pi}}\int_{-\infty}^a e^{\frac{-t^2}{2}}dt
\end{align*}
\end{thm}

\hf \hf The idea behind his proof is to first show that the conclusion of the theorem holds for a suitable class of step functions $\varphi$. He achieves this by showing that the random variables 
$X^{\varphi}_n$ are `statistically independent' (see \cite[Section 3]{Kac_Ber}). He then uses an approximation argument to show that the conclusion holds for a more general $f$ 
(whose Fourier coefficients satisfy the appropriate decay conditions). In particular, Kac establishes convergence in distribution to the standard 
normal distribution for a 
suitably normalized sample average $W^f_n$ of the collection $\{X^f_k\}_{k=0}^{\infty}$. \\ 
\hf \hf In this paper, we apply Stein's method to the question considered by Kac, namely to the collection of random variables 
$\{X^f_k\}_{k=0}^{\infty}$. We show that the same normalized average $W^f_n$ converges to the standard 
normal distribution in the \textbf{Wasserstein Metric}. 

\begin{defn}[{\cite[p. 214]{Rs_St},\cite{Villani-book}}]
\label{wass_metric}
Let $X, Y:\Omega\longrightarrow \mathbb{R}$ be two random variables. 
We define the \textbf{Wasserstein metric} between these two variables as 
\begin{align*}
d_W(X,Y)&:= \textnormal{sup}\{\mathbb{E}(h(X)-h(Y))|: h \in \textnormal{Lip}^{1}(\mathbb{R}, \mathbb{R})\}, 
\end{align*}
where $\textnormal{Lip}^{1}(\mathbb{R}, \mathbb{R})$ denotes the space of $1$-Lipschitz functions from $\mathbb{R}$ to $\mathbb{R}$.
\end{defn}
It is well-known that $d_W$ is indeed a metric and defines a notion of convergence, which further implies the convergence in 
distribution. We are now ready to state the main result of this paper. 
\begin{thm}[Main result: CLT for Angle Doubling Map]
\label{mt_bmp}
Let $f:[0,1]\longrightarrow \mathbb{R}$ be a non-zero $L^1$ function, with Fourier coefficients satisfying \eqref{Fourier_condition_f}. 
Let $\{X^f_k\}_{k=0}^{\infty}$ be the collection of random variables as defined by \eqref{Xf_defn} and let $W^f_n$ be as defined 
by \eqref{Wnf_defn}. Suppose equation \eqref{sigma_n_lim_assumption} is satisfied. Under these assumptions,   
the random variables $W^f_n$ converge to the standard normal random variable in the Wasserstein metric. 
\end{thm}
Since the convergence in Wasserstein metric is strictly stronger than the convergence in distribution (see \cite[p. 10]{Bobkov-LedouxBk}), we have an alternative proof of Kac's result. 

\section{A Recall of Preliminary Results and Proof of Theorem \ref{mt_bmp}}\label{S3:proofs}

\subsection{A Recall of Preliminary Results}
We review a few important concepts and results from \cite{Kac_Ber} and \cite{Rs_St}.

\begin{defn}[{\cite[Definition 3.5]{Rs_St}}]
Let $X_0, X_1, X_2, \ldots, X_n$ be a finite collection of random variables. We say this collection has dependency 
neighbourhoods $N_i \subset \{0, 1,2, \ldots, n\}$ for $i=0, 1, 2, \ldots n$ if $i \in N_i$ and $X_i$ is independent of $X_j$, whenever 
$j\notin N_i$.
\end{defn}

The central ingredient of our proof is to invoke the following result.

\begin{thm}[{\cite[Theorem 3.6, pp. 221]{Rs_St}}]
\label{ross}
Let $X_0, X_1, X_2, \ldots, X_n$ be random variables such that $\mathbb{E}[X_i^4]< \infty$, $\mathbb{E}[X_i] =0$. Let the collection $\{X_0, X_1, \ldots, X_n\}$ have dependency neighborhoods $N_i$, $i=0, \ldots, n$ 
and define $D_n:= \textnormal{max}_{0\leq i\leq n} |N_i|$. Then
\begin{align*}
d_W(W_n^X, Z) & \leq \frac{D_n^2}{\sigma_n^3} \sum_{i=1}^n \mathbb{E}|X_i|^3 + \frac{\sqrt{28} D_n^{\frac{3}{2}}}{\sqrt{\pi} \sigma_n^2}\sqrt{\sum_{i=1}^n\mathbb{E}|X_i|^4},
\end{align*}
where $\sigma_n^2$ is the variance of $X_0 + X_1 + X_2 + \ldots + X_n$, $W_n^X := \frac{X_0 + X_1 + X_2 + \ldots + X_n}{\sqrt{n} \sigma_n}$ and $Z$ denotes a standard normal random variable. 
\end{thm}
\hf \hf Let us now review the basic idea of Kac's proof. Kac showed that the main conclusion holds for a suitable class of step functions.

\begin{lmm}[{\cite[Lemma 5, pp. 40]{Kac_Ber}}]
\label{kac_step_function_thm}
Let $[0,1]$ be divided into $2^r$ intervals of equal length. Let $\varphi:[0,1]\longrightarrow \mathbb{R}$ be a step function such that 
on the interval $[\frac{i}{2^r}, \frac{i+1}{2^r}]$ it is a constant, for all $i=0$ to $2^r-1$. 
Assume $\varphi$ satisfies equation \eqref{sigma_n_lim_assumption}. 

\end{lmm} 
\begin{rem}
Note that we are not making any assumptions on the Fourier coefficients of $\varphi$. 
\end{rem}
\noindent Since this lemma will play a very important role in our paper, we briefly recall the proof from Kac's paper 
\cite{Kac_Ber}. Given a $t\in [0,1]$, 
let us consider the binary expansion of $t$. Let us consider the digits of $t$ to be $\varepsilon_1(t), \varepsilon_2(t), \ldots$. 
Note that each $\varepsilon_i(t)$ is either $0$ or $1$. Kac made the crucial observation that the step functions considered here can be represented as
\begin{align}
\varphi(t)&= P(\varepsilon_1(t), \varepsilon_2(t), \ldots, \varepsilon_r(t)), \label{phi_t_is_poly}
\end{align}
where $P$ is a \textit{polynomial} in $r$-variables. In other words, $\varphi(t)$ depends \textit{only} on the first $r$ 
digits of $t$ and furthermore, this dependence is through some explicit polynomial $P$ as in \eqref{phi_t_is_poly}. Using this fact, Kac then established the statistical independence of the random variables $X^{\varphi}_n$. Finally, he invokes a Lemma by Markoff \cite[pp. 302-309]{Mk}, 
that says that if the 
random variables are statistically independent, then the conclusion of the Central Limit Theorem still holds.  
This fact is also stated in \cite[Lemma 4, pp. 38-39]{Kac_Ber}, by referring to Markoff's paper \cite{Mk}. \\ 
\hf \hf After that, the author uses an approximation  
argument (namely he approximates $f$ by a step function) 
to show that Theorem \ref{mt_kac} is true for a general class of functions $f$ (provided of course, $f$ satisfies the hypothesis of the Theorem).

\subsection{Proof of Theorem \ref{mt_bmp}}
We first establish the result for a suitable step function, using the result on dependency neighborhood, namely Theorem \ref{ross}. This gives an alternative proof of Lemma \ref{kac_step_function_thm}. For clarity, we restate the result.

\begin{lmm}
\label{step_wass}
Let $[0,1]$ be divided into $2^r$ intervals of equal length. Let $\varphi:[0,1]\longrightarrow \mathbb{R}$ be a step function such that 
on the interval $[\frac{i}{2^r}, \frac{i+1}{2^r}]$ it is a constant, for all $i=0$ to $2^r-1$. 
Assume $\varphi$ satisfies equation \eqref{sigma_n_lim_assumption}. Then the conclusion of Theorem \ref{mt_bmp} holds. 
\end{lmm}
\begin{rem}
    Again, no assumptions are being made on the Fourier coefficients of $\varphi$.
\end{rem}
\noindent \textbf{Proof of Lemma \ref{step_wass}:} We will be using Theorem \ref{ross} to prove this Lemma. 
First of all, we note that without loss in generality, we can assume that 
\begin{align}
\int\limits_{0}^{1} \varphi(t) dt & = 0. \label{int_varphi_is_zero}
\end{align}
This is because, if the integral is equal to $\mu$, which is non-zero, then we can define the step function
\begin{align*}
\tilde{\varphi}(t)&:= \varphi(t)-\mu. 
\end{align*}
Assuming we have proved the Lemma for $\tilde{\varphi}$, the result follows for $\varphi$ as well. Henceforth, we will assume \eqref{int_varphi_is_zero}. 
Next, given a $t\in [0,1]$, the digits of $t$ in the binary expansion of $t$. By equation \eqref{phi_t_is_poly}, $\varphi(t)$ depends \textit{only} on the first $r$ digits of $t$, namely it depends 
only on $\varepsilon_1(t), \ldots, \varepsilon_r(t)$.\\ 
\hf \hf We now note that for all $t$,  $n$ and $k$, 
\begin{align}
\varepsilon_k(T^n t) & = \varepsilon_{k+n}(t). \label{ep_n_k} 
\end{align}
Furthermore, the random variables $\varepsilon_{i}:[0,1]\longrightarrow \mathbb{R}$ are i.i.d. and 
\begin{align*}
X^{\varphi}_n(t)&= \varphi(T^n t) \\ 
                & = P(\varepsilon_{1+n}(t), \varepsilon_{2+n}(t), \ldots, \varepsilon_{r+n}(t)). 
\end{align*}
Since, $\varepsilon_i$'s are i.i.d., the joint distribution of 
$(\varepsilon_{1 + n}, \cdots, \varepsilon_{r + n})$ depends on $r$ and hence on $\varphi$, but not on $n$. 
Hence $X_n^\varphi$'s are identically distributed. In particular, the variance of $X_n^\varphi$ is 
same as variance of $X_0^\varphi$, which is equal to the integral of $\varphi^2$ from $0$ to $1$, i.e. 
\begin{align*}
\textnormal{Var}(X_n^\varphi)&= \int\limits_{0}^1\varphi(t)^2 dt.  
\end{align*}
For the same reason, the third and fourth absolute moments, namely 
$\mathbb{E}|X_n^\varphi|^3$ and $\mathbb{E}|X_n^\varphi|^4$ are constants, independent of $n$.
We also note that since $\varphi$ is a step function; we have
\begin{align}
\int\limits_{0}^{1} \varphi(t)^4\, dt < \infty, \label{phi_fourth_moment}
\end{align}
which ensures the finiteness of $\mathbb{E}|X_n^\varphi|^4$. We have already assumed that 
the integral of $\varphi$ is zero (equation \eqref{int_varphi_is_zero}). This implies that 
all the expectations $\mathbb{E}[X_n^{\varphi}]$ are equal to zero 
(since they are all equal to $\mathbb{E}[X_0^{\varphi}]$, which in turn is equal to the integral of $\varphi$). 
Hence, the first two hypotheses of Theorem \ref{ross} are satisfied. \\ 
\hf\hf Let us now study the dependency neighborhood of the random variables $X^{\varphi}_n$. Let us consider the collection 
\[ \{X^{\varphi}_0, X^{\varphi}_1, \ldots, X^{\varphi}_n\}.  \]
Define $N_k$ to be the following set:
\begin{align*}
N_k&:= \{k-r+1, \ldots, k+r-1 \}.
\end{align*}
We now note that $X^{\varphi}_k$ is independent of the $X^{\varphi}_j$ if $j \notin N_k$. 
This is because the $\varepsilon_i$'s in $X^{\varphi}_k$ do not appear in the other random variable $X^{\varphi}_j$ if $j\notin N_k$. 
Moreover, $N_k$ is the largest such set of indices with this property. Therefore the dependency neighborhood of $X_k^\varphi$ is given by $N_k$. Since $|N_k|$ is atmost $2r - 1$ for all $k$, we have 
\begin{align}
D_n &\leq 2r - 1, \label{dep_nhbd}
\end{align}
a constant independent of $n$ and which depends only on $\varphi$ through $r$.\\ 
\hf \hf Let us now compute the covariance of $X^{\varphi}_i$ and $X^{\varphi}_j$.   
We claim that 
\begin{align}
\textnormal{Cov}(X_i^\varphi, X_j^\varphi) & = \rho(|i - j|)\int\limits_{0}^1 \varphi(t)^2 dt, \label{claim1}
\end{align}
where $\rho(i)$ denotes the correlation between $X_0^\varphi$ and $X_i^\varphi$. 
Let us denote by $\stackrel{d}{=}$ the equality of two random vectors in distribution. 
Without loss of generality, take $i < j$. We have, 
\begin{align}                         
(X_i^\varphi, X_j^\varphi) &= (P(\varepsilon_{1 + i}, \cdots, \varepsilon_{r + i}), P(\varepsilon_{1 + j}, \cdots, 
\varepsilon_{r + j})) \label{ab1}\\                                                                                                                                                                               
&\stackrel{d}{=} (P(\varepsilon_{1}, \cdots, \varepsilon_{r}), P(\varepsilon_{1 + j - i}, \cdots, \varepsilon_{r + j - i}))
\label{ab2}
\end{align}
which implies
\begin{equation}\label{ab3}
(X_i^\varphi, X_j^\varphi) \stackrel{d}{=} (X_0^\varphi, X_{j - i}^\varphi).
\end{equation}
To avoid confusion, let us clarify the meaning of the right hand sides of equations \eqref{ab1} and \eqref{ab2}. 
The expression $P(\varepsilon_{1}, \cdots, \varepsilon_{r})$ denotes the random variable $t$ going to 
$P(\varepsilon_{1}(t), \cdots, \varepsilon_{r}(t))$.\\  
\hf \hf Using equation \eqref{ab3}, we conclude that 
\begin{align*}
\textnormal{Cov}(X_i^\varphi, X_j^\varphi) & = 
\textnormal{Cov}(X_0^\varphi, X_{j - i}^\varphi) = \rho(j - i) 
\sqrt{\textnormal{Var}(X_0^\varphi) \textnormal{Var}(X_{j - i}^\varphi)},
\end{align*}
which proves \eqref{claim1}. \qed \\ 
\hf \hf Finally, we estimate the growth of the sample standard deviation, namely $\sigma_n$. Note that 
\begin{align}
\frac{1}{n}\sigma_n^2 &= \frac{1}{n} \left(\textnormal{Var}\Big(\sum\limits_{i=0}^{n-1} X_i^\varphi\Big)\right) \nonumber \\
&= \frac{1}{n } \left(\sum_{i = 0}^{n-1} \textnormal{Var}(X_i^\varphi) + \sum_{\stackrel{i, j \in \{0, 1,\ldots, n\}}{i \neq j}} \textnormal{Cov}(X_i^\varphi, X_j^\varphi)\right) \nonumber \\ 
& = \left(1 + \frac{J_n}{n} \right)\int\limits_{0}^{1} \varphi(x)^2\, dx, \qquad \textnormal{where} \qquad 
J_n:= \sum_{\stackrel{i, j \in \{0, 1,\ldots, n-1\}}{i \neq j}} \rho(|i - j|).  \label{cd1} 
\end{align}
Let us now analyse $J_n$. Whenever, $|i - j| > r - 1$, $X^{\varphi}_0$ and $X^{\varphi}_{|i-j|}$ are independent and consequently, $\rho(|i-j|) = 0$. Then,
\[J_n:= \sum_{i = 0}^{n-1} J_{n, i} \qquad \textnormal{where} \qquad J_{n, i} :=  \sum_{\stackrel{j: 0 < |i - j| \leq r - 1}{j \in \{0, 1,\ldots, n-1\}}} \rho(|i - j|).\]
Without loss of generality, we work with $n \geq 2r - 1$. For $i \geq r - 1$, $J_{n , i} = 2\sum\limits_{k=1}^{r-1} \rho(k) =: C_3$, a constant independent of $n$. Again, by construction, $|\rho(j)| \leq 1$ for any $j$ and consequently, $\sum_{i = 0}^{r-2} J_{n, i}$ is bounded by $(r-1)^2$, independent of $n$.
Hence, 
\[\lim_{n \to \infty} \frac{J_n}{n} = \lim_{n \to \infty} \frac{1}{n} \sum_{i = 0}^{r-2} J_{n, i} + \lim_{n \to \infty} \frac{1}{n} \sum_{i = r - 1}^{n - 1} J_{n, i} = 0 + C_3 = C_3.\]
Using equation \eqref{cd1}, we have 
\[\lim_{n\to \infty}\frac{\sigma_n^2}{n} = (1+C_3) \int_0^1 \varphi(x)^2\, dx. \]

We now conclude, using Theorem \ref{ross} that $d_W(W_n^{X_{\varphi}}, Z) = O(\frac{1}{\sqrt{n}})$ as $n$ goes to infinity. 
This completes the proof of Lemma \ref{step_wass}. \qed \\
\hf \hf We are now ready to prove our main Theorem for a general function. First, let us prove one auxiliary result. 
We make a very simple observation that bounds the size of Wasserstein distance by the size of the $L^2$ distance.
\begin{lmm}
\label{wass_vs_l2}
Let $X$ and $Y$ be random variables. Then 
\begin{align*}
d_W(X, Y) & \leq ||X-Y||_{L^2}. 
\end{align*}
\end{lmm}
\noindent \textbf{Proof:} The proof follows by unwinding definitions. We note that 
\begin{align*}
d_W(X, Y) & = \textnormal{sup}\{ |\mathbb{E}(h(X)-h(Y))|: ~~\textnormal{$h$ is $1$-Lipschitz} \} \\ 
          & \leq \textnormal{sup}\{ \mathbb{E}|(h(X)-h(Y)|): ~~\textnormal{$h$ is $1$-Lipschitz} \} \\ 
          & \leq \mathbb{E}|X-Y| \qquad \textnormal{using the fact that the $h$ is $1$-Lipschitz}\\ 
          &\leq \Big(\mathbb{E}|X-Y|^2 \Big)^{\frac{1}{2}} = ||X-Y||_{L^2}. \qquad \qed 
\end{align*}
\noindent We are now ready to prove the main theorem. \\

\noindent \textbf{Proof of Theorem \ref{mt_bmp}:} For the function $f$ under consideration, we recall from \cite[Lemma 3]{Kac_Ber} that given an $\varepsilon>0$, there exists a step function $\varphi$ with its such that the endpoints of its intervals of constancy of the form $\frac{i}{2^r}$, $\int_0^1 \varphi(x) \, dx = 0$ and 
\[\lim_{m \to \infty} \frac{1}{m} \left\|\sum_{k = 0}^{m-1} \varphi(T^k\cdot) \right\|^2_{L^2} > 0\]
such that for large $n$,
\begin{align}
||W_n^f-W_n^{\varphi}||_{L^2} &< \frac{\varepsilon}{2}. \label{eq1}
\end{align}
By the triangle inequality for the Wasserstein metric $d_W$ and using Lemmas \ref{step_wass} and \ref{wass_vs_l2}, we have
\[d_W(W_n^f, Z) \leq d_W(W_n^f, W_n^{\varphi}) + d_W(W_n^{\varphi}, Z) \leq ||W_n^f-W_n^{\varphi}||_{L^2} + d_W(W_n^{\varphi}, Z) < \varepsilon\]
for large $n$. This concludes the proof. \qed

\section{Acknowledgement} 
We gratefully acknowledge constructive comments from an anonymous reviewer. Suprio Bhar gratefully acknowledges support from the DST INSPIRE Faculty Scheme (DST/INSPIRE/04/2017/002835) from Government of India.

\noindent\textbf{Declarations of interest:} none.




\end{document}